\documentclass[12pt]{article}
\usepackage{amsfonts}
\usepackage{a4wide}
\usepackage{amsmath,amscd,amsthm,a4,amssymb}
\usepackage{sectsty}

\sectionfont{\centering}
\numberwithin{equation}{section}

\begin{document}

\begin{center}
{\LARGE \textbf{On the Kurepa and inhomogeneous Cauchy functional equations}}

\

\textbf{Rashid A. Aliev}$^{1}$ \textbf{and Vugar E. Ismailov}\footnote{%
Corresponding author}$^{2}$ \vspace{3mm}

$^{1,2}${Institute of Mathematics and Mechanics, Baku, Azerbaijan} \vspace{%
1mm}

$^{1}${Baku State University, Baku, Azerbaijan} \vspace{1mm}

$^{1,2}${Center for Mathematics and its Applications, Khazar University,
Baku, Azerbaijan} \vspace{1mm}

e-mail: $^{1}$aliyevrashid@mail.ru, $^{2}$vugaris@mail.ru
\end{center}

\bigskip

\textbf{Abstract.} It follows from de Bruijn's results that if a continuous
or $k$-th order continuously differentiable function $F(x,y)$ is a solution
of the Kurepa functional equation, then it can be expressed as $%
F(x,y)=f(x+y)-f(x)-f(y)$ with the continuous $f$ or the $k$-th order
continuously differentiable $f$, respectively. These two facts strengthen
the corresponding results of Kurepa and Erd\"{o}s. In this paper, we provide
new and constructive proofs for these facts. In addition to practically
useful recipes given here for construction of $f$, we also estimate its
modulus of continuity.

\bigskip


\smallskip


\bigskip

\section{Introduction}

There are extensive investigations on the Kurepa functional equation
\begin{equation}
F(x+y,z)+F(x,y)=F(y,z)+F(x,y+z)  \label{1.1}
\end{equation}%
with the aim to find the general solution of this equation under certain
conditions. If we consider the solutions of (1.1) as a functional equation
then we get substantially the inhomogeneous Cauchy functional equation
\begin{equation}
F(x,y)=g(x+y)-g(x)-g(y).  \label{1.2}
\end{equation}%
Clearly, any function of the form \eqref{1.2} satisfies \eqref{1.1} and it
was first proved by Kurepa \cite{K} that all differentiable solutions of %
\eqref{1.1} are of the form \eqref{1.2}. Erd\"{o}s \cite{Er} showed that the
last assertion is valid also for continuous solutions of \eqref{1.1}. The
question whether \eqref{1.2} is the general solution of \eqref{1.1} has a
negative answer (see \cite{Er}).

The equation \eqref{1.2} has drawn the attention of many authors and been
studied for various spaces and forms of $F$ (see e.g., \cite%
{Bor,Br2,D,E3,Er,Fen,Jr,Jes}). There are several natural and useful
generalizations of the inhomogeneous Cauchy functional equation (see e.g.
\cite{Br,E,E2,Pr}).

Note that the above results of Kurepa and Erd\"{o}s do not say anything
about the nature of $g$ in \eqref{1.2}. However, it follows from de Bruijn's
remarkable results that if $F(x,y)$ is a solution of \eqref{1.1} from the
class $C^{k}(\mathbb{R}^{2})$ or $C(\mathbb{R}^{2})$, then it can be written
in the form \eqref{1.2} with $g$ from $C^{k}(\mathbb{R})$ or $C(\mathbb{R})$%
, respectively. To see this, recall that for a function $g:\mathbb{R}%
\rightarrow \mathbb{R}$ and any $h\in \mathbb{R}$ the difference function $%
\Delta _{h}g:\mathbb{R}\rightarrow \mathbb{R}$ is defined as
\begin{equation*}
\Delta _{h}g(x)=g(x+h)-g(x).
\end{equation*}%
Let $\mathcal{F}$ be a class of functions defined on $\mathbb{R}$. The class
$\mathcal{F}$ is said to have the difference property if every function $g:%
\mathbb{R}\rightarrow \mathbb{R}$, for which $\Delta _{h}g\in \mathcal{F}$
for each $h\in \mathbb{R}$, is of the form $g=f+A$, where $f\in \mathcal{F}$
and $A$ is an additive function (see \cite{1}). A function $A$ is called
additive if it satisfies the Cauchy functional equation $A(x+y)=A(x)+A(y)$.
In 1951, de Bruijn \cite{1} showed that the class of continuous functions on
$\mathbb{R}$ has the difference property, hence resolving Erd\"{o}s's famous
conjecture. He also proved that the difference property holds for various
essential function classes (see \cite{1,2}), in particular, for the class $%
C^{k}(\mathbb{R)}$, functions with continuous derivatives up to order $k$.
Now if $F$ is of the form (1.2) and $F\in C^{k}(\mathbb{R}^{2})$, then $%
g(x+y)-g(x)$, as a function of $x$, belongs to $C^{k}(\mathbb{R})$ for every
fixed $y$. Since $C^{k}(\mathbb{R)}$ has the difference property, we can
write that $g=f+A,$ where $f\in C^{k}(\mathbb{R)}$ and $A$ is additive.
Hence $F(x,y)=g(x+y)-g(x)-g(y)=f(x+y)-f(x)-f(y)$. Using the difference
property of $C(\mathbb{R)}$, it can be shown by the same way that if $F$ is
of the form (1.2) and $F\in C(\mathbb{R}^{2})$, then there exists $f\in C(%
\mathbb{R})$ such that $F(x,y)=f(x+y)-f(x)-f(y)$. Note that these results
establish only existence of such $f$. We know nothing about intrinsic nature
and properties of $f$. In this paper, we provide different and constructive
proofs for these results. Our method allows not only to construct $f$ in
each considered case, but also to estimate its modulus of continuity.

\bigskip

\section{Solutions of (1.1) from the class $C^{k}(\mathbb{R}^{2})$}

In this section, we see that if $F$ is a solution of the Kurepa functional
equation \eqref{1.1} from the class $C^{k}(\mathbb{R}^{2}),$ then it can be
written in the form \eqref{1.2} with $g\in C^{k}(\mathbb{R})$. First we
prove the following theorem, which deals with the inhomogeneous Cauchy
functional equation.

\bigskip

\textbf{Theorem 2.1.} \textit{Assume that a function $F\in C^{k}(\mathbb{R}%
^{2})$ has the form}

\begin{equation}
F(x,y)=g(x+y)-g(x)-g(y),  \label{2.0}
\end{equation}%
\textit{where $g$ is an arbitrarily behaved univariate function and $k\geq
1. $ Then $F$ can be represented also in the form}
\begin{equation*}
F(x,y)=f(x+y)-f(x)-f(y),
\end{equation*}%
\textit{where the function $f\in C^{k}(\mathbb{R})$.}

\bigskip

\textbf{Proof.} By $\Delta _{l}^{(\delta )}f$ we denote the increment of a
function $f$ in a direction $l=(l^{\prime },l^{\prime \prime }).$ That is,
\begin{equation*}
\Delta _{l}^{(\delta )}f(x,y)=f(x+l^{\prime }\delta ,y+l^{\prime \prime
}\delta )-f(x,y).
\end{equation*}%
We also use the notation $\frac{\partial f}{\partial l}$ which denotes the
derivative of $f$ in the direction $l$.

It is easy to check that the increment of a function $f(x,y)$ of the form $%
g(ax+by)$ in a direction perpendicular to $(a,b)$ is zero. Let $l=(1/\sqrt{2}%
,-1/\sqrt{2})$, $e_{1}=(1,0)$ and $e_{2}=(0,1)$. Then for any number $\delta
\in \mathbb{R}$, we have

\begin{equation}
\Delta _{l}^{(\delta )}F(x,y)=\Delta _{l}^{(\delta )}\left[ -g(x)-g(y)\right]
.  \label{2.3}
\end{equation}

Denote the left hand side of \eqref{2.3} by $S(x,y).$ That is, set%
\begin{equation*}
S(x,y)\overset{def}{=}\Delta _{l}^{(\delta )}F(x,y).
\end{equation*}%
Then from \eqref{2.3} it follows that for any real numbers $\delta _{1}$and $%
\delta _{2}$,

\begin{equation*}
\Delta _{e_{1}}^{(\delta _{1})}\Delta _{e_{2}}^{(\delta _{2})}S(x,y)=0,
\end{equation*}%
or in expanded form,

\begin{equation*}
S(x+\delta _{1},y+\delta _{2})-S(x,y+\delta _{2})-S(x+\delta
_{1},y)+S(x,y)=0.
\end{equation*}%
Putting in the last equality $\delta _{1}=-x,$ $\delta _{2}=-y$, we obtain
that

\begin{equation*}
S(x,y)=S(x,0)+S(0,y)-S(0,0).
\end{equation*}%
This means that

\begin{equation*}
\Delta _{l}^{(\delta )}F(x,y)=\Delta _{l}^{(\delta )}F(x,0)+\Delta
_{l}^{(\delta )}F(0,y)-\Delta _{l}^{(\delta )}F(0,0).
\end{equation*}

By the hypothesis of the theorem, the directional derivative $\frac{\partial
F}{\partial l}$ exists at any point $(x,y)\in $ $\mathbb{R}^{2}$. Thus, it
follows from the above formula that

\begin{equation}
\frac{\partial F}{\partial l}(x,y)=h_{1}(x)+h_{2}(y),  \label{2.4}
\end{equation}%
where $h_{1}(x)=\frac{\partial F}{\partial l}(x,0)$ and $h_{2}(y)=\frac{%
\partial F}{\partial l}(0,y)-\frac{\partial F}{\partial l}(0,0)$. Note that $%
h_{1}$ and $h_{2}$ belong to the class $C^{k-1}(\mathbb{R}).$

It follows from \eqref{2.0} that $F(x,y)$ is a symmetric function (that is, $%
F(x,y)=F(y,x)$ holds for any $x,y$). Therefore, $\frac{\partial F}{\partial l%
}(t,0)=-\frac{\partial F}{\partial l}(0,t)$ for any $t$. On the other hand,
it is easy to verify that

\begin{equation*}
F(\delta ,0)-F(0,0)=F(0,\delta )-F(0,0)=0.
\end{equation*}%
Consequently,

\begin{equation*}
\frac{\partial F}{\partial x}(0,0)=\frac{\partial F}{\partial y}(0,0)=0
\end{equation*}%
and hence $\frac{\partial F}{\partial l}(0,0)=0$. Thus, we see that $%
h_{1}=-h_{2}.$

\bigskip

By $u_{1}$ and $u_{2}$ denote the antiderivatives of $h_{1}$ and $h_{2}$
satisfying the conditions $h_{1}(0)=h_{2}(0)=0$ and multiplied by the
numbers $1/(e_{1}\cdot l)=\sqrt{2}$ and $1/(e_{2}\cdot l)=-\sqrt{2}$
correspondingly (here, $e\cdot l$ denotes the scalar product between the
vectors $e$ and $l$). That is,
\begin{eqnarray*}
u_{1}(x) &=&\sqrt{2}\int_{0}^{x}h_{1}(z)dz; \\
u_{2}(y) &=&-\sqrt{2}\int_{0}^{y}h_{2}(z)dz.
\end{eqnarray*}%
Obviously, the function

\begin{equation*}
F_{1}(x,y)=u_{1}(x)+u_{2}(y)
\end{equation*}%
obeys the equality

\begin{equation}
\frac{\partial F_{1}}{\partial l}(x,y)=h_{1}(x)+h_{2}(y).  \label{2.5}
\end{equation}%
From \eqref{2.4} and \eqref{2.5} we obtain that

\begin{equation*}
\frac{\partial }{\partial l}\left[ F-F_{1}\right] =0.
\end{equation*}%
Hence, for some function $\varphi (x+y),$

\begin{equation}
F(x,y)=F_{1}(x,y)+\varphi (x+y)=u_{1}(x)+u_{2}(y)+\varphi (x+y)  \label{2.6}
\end{equation}%
Here all the functions $u_{1},u_{2},\varphi \in C^{k}(\mathbb{R}).$ Since $%
h_{1}=-h_{2}$, it follows that $u_{1}=u_{2}.$ Set $u=-u_{1}$. Then %
\eqref{2.6} can be written in the form

\begin{equation}
F(x,y)=\varphi (x+y)-u(x)-u(y).  \label{2.7}
\end{equation}

Let us show that $F$ can also be written in the form

\begin{equation}
F(x,y)=f(x+y)-f(x)-f(y)  \label{2.8}
\end{equation}%
with $f\in C^{k}(\mathbb{R}).$ Indeed, taking $y=0$, we obtain from %
\eqref{2.0} and \eqref{2.7} that

\begin{equation*}
-g(0)=\varphi (x)-u(x)-u(0).
\end{equation*}%
Hence,
\begin{equation*}
\varphi (x)=u(x)+C,
\end{equation*}%
where the constant $C=u(0)-g(0)$. Consider now the function $f=u-C$. We see
that \eqref{2.8} holds for this function. The theorem has been proved.

\bigskip

Now we are ready to formulate the second result of this section, which
follows from Theorem 2.1 and the theorem of Kurepa that all differentiable
solutions of \eqref{1.1} are of the form \eqref{1.2}.

\bigskip

\textbf{Theorem 2.2.} \textit{If $F(x,y)$ is a solution of \eqref{1.1} from
the class $C^{k}(\mathbb{R}^{2})$, $k\geq 1$, then it can be written in the
form \eqref{1.2} with $g\in C^{k}(\mathbb{R})$.}

\bigskip

\section{Continuous solutions of (1.1)}

Assume $f$ is a function given on a set $E\subset \mathbb{R}^{d}$, $d\geq 1$%
. The modulus of continuity of $f$ on $E$ is defined as%
\begin{equation*}
\omega (f;\delta ;E)=\sup \left\{ \left\vert f(x)-f(y)\right\vert :x,y\in E,%
\text{ }\left\vert x-y\right\vert \leq \delta \right\} ,\text{ }0\leq \delta
\leq diamE.
\end{equation*}

In this section we show that if $F$ is a continuous solution of the Kurepa
functional equation \eqref{1.1}, then it can be written in the form %
\eqref{1.2} with the continuous $g$. First we prove the following theorem,
which establishes the existence and gives a recipe for construction of a
continuous solution of the inhomogeneous Cauchy functional equation. Our
methods also allows to estimate its modulus of continuity.

\bigskip

\textbf{Theorem 3.1.} \textit{Assume a function $F\in C(\mathbb{R}^{2})$ has
the form}
\begin{equation}
F(x,y)=g(x+y)-g(x)-g(y),  \label{3.1}
\end{equation}%
\textit{where $g$ is an arbitrarily behaved function. Then there exists a
function $f\in C(\mathbb{R})$ such that}
\begin{equation}
F(x,y)=f(x+y)-f(x)-f(y).  \label{3.0}
\end{equation}%
\textit{For each real $t$ the function $f(t)$ can be obtained from $g$ by
taking a limit of $g(t_{n}),$ as $t_{n}\rightarrow t,$ $t_{n}\in \mathbb{Q}$%
. In addition the following inequality holds}
\begin{equation}
\omega \left( f;\delta ;[-M,M]\right) \leq 3\omega \left( F;\delta
;[-M,M]^{2}\right),  \label{3.20}
\end{equation}%
\textit{where $\delta \in \left( 0,\frac{1}{2}\right) \cap \mathbb{Q}$ and $%
M\geq 1$.}

\bigskip

\textbf{Proof.} Write Eq. \eqref{3.1} in the form
\begin{equation}
H(x,y)=h(x+y)-h(x)-h(y),  \label{3.3}
\end{equation}%
where
\begin{equation*}
H(x,y)=F(x,y)+g(0),
\end{equation*}%
and
\begin{equation*}
h(t)=g(t)-g(0).
\end{equation*}%
Since $h(0)=0,$ we obtain from \eqref{3.3} that

\begin{equation}
H(x,0)=H(0,y)=0.  \label{3.4}
\end{equation}%
It is not difficult to verify that for any nonnegative integer $k$%
\begin{equation*}
H(x,(k-1)x)=h(kx)-h(x)-h((k-1)x).
\end{equation*}%
It follows from this formula that

\begin{equation*}
h(kx)=kh(x)+\sum_{i=1}^{k-1}H(x,ix).
\end{equation*}%
Therefore, for any nonnegative integer $k$ we have

\begin{equation}
h(x)=\frac{1}{k}h(kx)-\frac{1}{k}\sum_{i=1}^{k-1}H(x,ix).  \label{3.5}
\end{equation}

Assume $\frac{p}{n}$ is any simple fraction from the interval $(0,\frac{1}{2}%
)$. Put $m_{0}=\left[ \frac{n}{p}\right] $, the whole number part of $\frac{n%
}{p}$. Obviously, $m_{0}\geq 2$ and the remainder $p_{1}=n-m_{0}p<p.$
Putting $x=\frac{p}{n}$ and $k=m_{0}$ in \eqref{3.5} yields

\begin{equation}
h\left( \frac{p}{n}\right) =\frac{1}{m_{0}}h\left( 1-\frac{p_{1}}{n}\right) -%
\frac{1}{m_{0}}\sum_{i=1}^{m_{0}-1}H\left( \frac{p}{n},\frac{ip}{n}\right) .
\label{3.6}
\end{equation}%
On the other hand, since

\begin{equation*}
H\left( \frac{p_{1}}{n},1-\frac{p_{1}}{n}\right) =h(1)-h\left( \frac{p_{1}}{n%
}\right) -h\left( 1-\frac{p_{1}}{n}\right) ,
\end{equation*}%
it follows from \eqref{3.6} that

\begin{equation}
h\left( \frac{p}{n}\right) =\frac{h(1)}{m_{0}}-\frac{1}{m_{0}}\left[
\sum_{i=1}^{m_{0}-1}H\left( \frac{p}{n},\frac{ip}{n}\right) +H\left( \frac{%
p_{1}}{n},1-\frac{p_{1}}{n}\right) \right] -\frac{1}{m_{0}}h\left( \frac{%
p_{1}}{n}\right) .  \label{3.7}
\end{equation}

Put now $m_{1}=\left[ \frac{n}{p_{1}}\right] $, $p_{2}=n-m_{1}p_{1}.$
Clearly, $0\leq p_{2}<p_{1}$. Similar to \eqref{3.7}, we can write that

\begin{equation*}
h\left( \frac{p_{1}}{n}\right) =\frac{h(1)}{m_{1}}-\frac{1}{m_{1}}\left[
\sum_{i=1}^{m_{1}-1}H\left( \frac{p_{1}}{n},\frac{ip_{1}}{n}\right) +H\left(
\frac{p_{2}}{n},1-\frac{p_{2}}{n}\right) \right] -\frac{1}{m_{1}}h\left(
\frac{p_{2}}{n}\right).
\end{equation*}

We can continue this process by defining the chain of pairs $(m_{2},p_{3}),$
$(m_{3},p_{4})$ until the pair $(m_{k-1},p_{k})$ with $p_{k}=0$ and writing
out the corresponding formulas for each pair. For example, the last formula
(that is, the formula corresponding to $(m_{k-1},p_{k})$) will be of the
form
\begin{equation*}
h\left( \frac{p_{k-1}}{n}\right) =\frac{h(1)}{m_{k-1}}
\end{equation*}%
\begin{equation}
-\frac{1}{m_{k-1}}\left[ \sum_{i=1}^{m_{k-1}-1}H\left( \frac{p_{k-1}}{n},%
\frac{ip_{k-1}}{n}\right) +H\left( \frac{p_{k}}{n},1-\frac{p_{k}}{n}\right) %
\right] -\frac{1}{m_{k-1}}h\left( \frac{p_{k}}{n}\right) .  \label{3.9}
\end{equation}%
Considering \eqref{3.9} in the formula corresponding to $(m_{k-2},p_{k-1})$,
then the obtained formula in the formula corresponding to $(m_{k-3},p_{k-2})$
and so on, we will arrive at the equality

\begin{equation*}
h\left( \frac{p}{n}\right) =h(1)\left[ \frac{1}{m_{0}}-\frac{1}{m_{0}m_{1}}%
+\cdot \cdot \cdot +\frac{(-1)^{k-1}}{m_{0}m_{1}\cdot \cdot \cdot m_{k-1}}%
\right]
\end{equation*}

\begin{equation*}
-\frac{1}{m_{0}}\left[ \sum_{i=1}^{m_{0}-1}H\left( \frac{p}{n},\frac{ip}{n}%
\right) +H\left( \frac{p_{1}}{n},1-\frac{p_{1}}{n}\right) \right]
\end{equation*}

\begin{equation*}
+\frac{1}{m_{0}m_{1}}\left[ \sum_{i=1}^{m_{1}-1}H\left( \frac{p_{1}}{n},%
\frac{ip_{1}}{n}\right) +H\left( \frac{p_{2}}{n},1-\frac{p_{2}}{n}\right) %
\right]
\end{equation*}

\begin{equation}
+\cdot \cdot \cdot +\frac{(-1)^{k}}{m_{0}m_{1}\cdot \cdot \cdot m_{k-1}}%
\sum_{i=1}^{m_{k-1}-1}H\left( \frac{p_{k-1}}{n},\frac{ip_{k-1}}{n}\right) .
\label{3.10}
\end{equation}%
Taking into account \eqref{3.4}, we obtain from \eqref{3.10} that

\begin{equation*}
\left\vert h\left( \frac{p}{n}\right) \right\vert \leq \left[ \frac{1}{m_{0}}%
-\frac{1}{m_{0}m_{1}}+\cdot \cdot \cdot +\frac{(-1)^{k-1}}{m_{0}m_{1}\cdot
\cdot \cdot m_{k-1}}\right] \left\vert h(1)\right\vert
\end{equation*}

\begin{equation}
+\left[ 1+\frac{1}{m_{0}}+\cdot \cdot \cdot +\frac{1}{m_{0}\cdot \cdot \cdot
m_{k-2}}\right] \omega \left( H;\frac{p}{n};[0,1]^{2}\right) .  \label{3.11}
\end{equation}%
Since $m_{0}\leq m_{1}\leq \cdot \cdot \cdot \leq m_{k-1},$ it is easy to
see that

\begin{equation*}
\frac{1}{m_{0}}-\frac{1}{m_{0}m_{1}}+\cdot \cdot \cdot +\frac{(-1)^{k-1}}{%
m_{0}m_{1}\cdot \cdot \cdot m_{k-1}}\leq \frac{1}{m_{0}}
\end{equation*}%
and

\begin{equation*}
1+\frac{1}{m_{0}}+\cdot \cdot \cdot +\frac{1}{m_{0}\cdot \cdot \cdot m_{k-2}}%
\leq \frac{m_{0}}{m_{0}-1}.
\end{equation*}%
Considering the above two inequalities in \eqref{3.11} we obtain that

\begin{equation}
\left\vert h\left( \frac{p}{n}\right) \right\vert \leq \frac{\left\vert
h(1)\right\vert }{m_{0}}+\frac{m_{0}}{m_{0}-1}\omega \left( H;\frac{p}{n}%
;[0,1]^{2}\right) .  \label{3.12}
\end{equation}%
Since $m_{0}=\left[ \frac{n}{p}\right] \geq 2,$ it follows from \eqref{3.12}
that

\begin{equation}
\left\vert h\left( \frac{p}{n}\right) \right\vert \leq \frac{2p\left\vert
h(1)\right\vert }{n}+2\omega \left( H;\frac{p}{n};[0,1]^{2}\right) .
\label{3.13}
\end{equation}

Assume now $\delta \in \left( 0,\frac{1}{2}\right) \cap \mathbb{Q}$, $M\geq
1 $ and $x,x+\delta \in \left[ -M,M\right] \cap \mathbb{Q}.$ Using %
\eqref{3.13} we can write that

\begin{equation}
\left\vert h(x+\delta )-h(x)\right\vert \leq \left\vert h(\delta
)\right\vert +\left\vert H(x,\delta )\right\vert \leq 2\delta \left\vert
h(1)\right\vert +3\omega \left( H;\delta ;[-M,M]^{2}\right) .  \label{3.14}
\end{equation}%
From \eqref{3.14} and the definitions of $h$ and $H$ it follows that

\begin{equation}
\omega (g;\delta ;[-M,M]\cap \mathbb{Q})\leq 2\delta \left\vert
g(1)-g(0)\right\vert +3\omega \left( F;\delta ;[-M,M]^{2}\right) ,
\label{3.2}
\end{equation}%
for any $\delta \in \left( 0,\frac{1}{2}\right) \cap \mathbb{Q}$ and $M\geq
1 $.

Introduce now the function
\begin{equation*}
u(t)=g(t)-\left[ g(1)-g(0)\right] t.
\end{equation*}%
For this function $u(1)=u(0)$ and

\begin{equation}
F(x,y)=u(x+y)-u(x)-u(y).  \label{3.17}
\end{equation}%
By \eqref{3.2}, which also holds for $u$, the restriction of $u$ to $\mathbb{%
Q}$ is uniformly continuous on every interval $[-M,M]\cap \mathbb{Q}$. Let
us denote this restriction by $v$.

Assume $y$ is any real number and $\{y_{n}\}_{n=1}^{\infty }$ is any
sequence of rational numbers tending to $y$ and $M>0$ so that $y_{n}\in
\lbrack -M,M]$ for any $n\in \mathbb{N}$. Since $v$ is uniformly continuous
on $[-M,M]\cap \mathbb{Q}$, the sequence $\{v(y_{n})\}_{n=1}^{\infty }$ is a
Cauchy sequence. Hence a finite $\lim_{n\rightarrow \infty }v(y_{n})$
exists. It is not difficult to understand that this limit is independent of
the choice of $\{y_{n}\}_{n=1}^{\infty }$.

By $f$ denote the following extension of $v$ to the set of real numbers.
\begin{equation*}
f(y)=\left\{
\begin{array}{c}
v(y),\text{ if }y\in \mathbb{Q}\text{;} \\
\lim_{n\rightarrow \infty }v(y_{n}),\text{ if }y\in \mathbb{R}\backslash
\mathbb{Q}\text{ and }\{y_{n}\}\text{ is a sequence in }\mathbb{Q}\text{
tending to }y.%
\end{array}%
\right.
\end{equation*}%
The above arguments show that $f$ is well defined on $\mathbb{R}$.

Let us now prove that \eqref{3.0} holds and $f\in C(\mathbb{R})$. Take any
point $(x,y)\in \mathbb{R}^{2}$ and a sequence of points $%
\{(x_{n},y_{n})\}_{n=1}^{\infty }$ with rational coordinates converging to $%
(x,y)$. Since $v$ is the restriction of $u$ to $\mathbb{Q}$, we have from %
\eqref{3.17} that

\begin{equation}
F(x_{n},y_{n})=v(x_{n}+y_{n})-v(x_{n})-v(y_{n}),\text{ for all }n=1,2,...
\label{3.18}
\end{equation}%
Tending $n\rightarrow \infty $ in both sides of \eqref{3.18} we obtain that

\begin{equation*}
F(x,y)=f(x+y)-f(x)-f(y).
\end{equation*}

It remains to prove the continuity of $f$. Since $v(1)=v(0)$ it follows from %
\eqref{3.17} and \eqref{3.2} that for any $\delta \in \left( 0,\frac{1}{2}%
\right) \cap \mathbb{Q}$, $M\geq 1$ and any numbers $r_{1},r_{2}\in \lbrack
-M,M]\cap \mathbb{Q}$, $\left\vert r_{1}-r_{2}\right\vert \leq \delta ,$ the
following inequality holds

\begin{equation}
\left\vert v(r_{1})-v(r_{2})\right\vert \leq 3\omega \left( F;\delta
;[-M,M]^{2}\right) .  \label{3.19}
\end{equation}%
Consider now any real numbers $a$ and $b$ satisfying $a,b\in \lbrack -M,M]$,
$\left\vert a-b\right\vert \leq \delta $ and take two sequences $%
\{a_{n}\}_{n=1}^{\infty }\subset \lbrack -M,M]\cap \mathbb{Q}$, $%
\{b_{n}\}_{n=1}^{\infty }\subset \lbrack -M,M]\cap \mathbb{Q}$ with the
property $\left\vert a_{n}-b_{n}\right\vert \leq \delta ,$ $n=1,2,...,$ and
tending to $a$ and $b$, respectively. By \eqref{3.19},
\begin{equation*}
\left\vert v(a_{n})-v(b_{n})\right\vert \leq 3\omega \left( F;\delta
;[-M,M]^{2}\right) .
\end{equation*}%
Taking limits\ on both sides of the above inequality gives

\begin{equation*}
\left\vert f(a)-f(b)\right\vert \leq 3\omega \left( F;\delta
;[-M,M]^{2}\right) .
\end{equation*}%
The last inequality proves \eqref{3.20}. This means that $f$ is uniformly
continuous on $[-M,M]$ and hence continuous on the whole real line. Theorem
3.1 has been proved.

\bigskip

Now we formulate the second result of this section, which follows from
Theorem 3.1 and the theorem of Erd\"{o}s that all continuous solutions of %
\eqref{1.1} are of the form \eqref{1.2}.

\bigskip

\textbf{Theorem 3.2.} \textit{If $F(x,y)$ is a continuous solution of %
\eqref{1.1}, then it can be written in the form \eqref{1.2} with the
continuous $g$.}

\bigskip

\textbf{Remark.} Some techniques used in this paper were also implemented in
\cite{A} to prove the double difference property for the class of locally H%
\"{o}lder continuous functions.

\end{document}